\documentclass[11pt,a4paper]{article}

\usepackage{theorem,enumerate}
\usepackage{amsmath,latexsym,amssymb,amsfonts}
\usepackage{eucal}
\usepackage{color}
\usepackage{comment}
\usepackage{cite}

\theorembodyfont{\normalfont\slshape}
\newtheorem{thm}{Theorem}[section]

\newtheorem{prop}[thm]{Proposition}
\newtheorem{lem}[thm]{Lemma}

\newtheorem{fact}[thm]{Fact}

\newtheorem{claim}{Claim}[section]

\newtheorem{problem}{Problem}

\newcommand{\proof}{\medbreak\noindent\textit{Proof.}\quad}

\newcommand{\qed}{{$\quad\square$\vs{3.6}}}

\newcommand{\vs}[1]{\vspace*{#1 mm}}

\def\AA{{ \mathcal{A}}}

\def\HH{{ \mathcal{H}}}

\def\PP{{ \mathcal{P}}}

\bfseries\normalfont

\title{Ramsey-type results for path covers and path partitions}

\author{
Shuya Chiba$^{1}$\footnote{\texttt{e-mail:schiba@kumamoto-u.ac.jp}}\and \
Michitaka Furuya$^{2}$\footnote{\texttt{e-mail:michitaka.furuya@gmail.com}}\vs{5}\\
$^{1}$\textsl{Applied Mathematics, Faculty of Advanced Science and Technology,}\\
\textsl{Kumamoto University,}\\
\textsl{2-39-1 Kurokami, Kumamoto 860-8555, Japan}\\
$^{2}$\textsl{College of Liberal Arts and Sciences,}\\
\textsl{Kitasato University,}\\
\textsl{1-15-1 Kitasato, Minami-ku, Sagamihara, Kanagawa 252-0373, Japan}
}

\date{}

\begin{document}

\maketitle

\begin{abstract}
A family $\mathcal{P}$ of subgraphs of $G$ is called a {\it path cover} (resp. a {\it path partition}) of $G$ if $\bigcup _{P\in \mathcal{P}}V(P)=V(G)$ (resp. $\dot\bigcup _{P\in \mathcal{P}}V(P)=V(G)$) and every element of $\mathcal{P}$ is a path.
The minimum cardinality of a path cover (resp. a path partition) of $G$ is denoted by  ${\rm pc}(G)$ (resp. ${\rm pp}(G)$).
In this paper, we characterize the forbidden subgraph conditions assuring us that ${\rm pc}(G)$ (or ${\rm pp}(G)$) is bounded by a constant.
Our main results introduce a new Ramsey-type problem.
\end{abstract}

\noindent
{\it Key words and phrases.}
path cover number, path partition number, forbidden subgraph, Ramsey number

\noindent
{\it AMS 2020 Mathematics Subject Classification.}
05C38, 05C55.

\section{Introduction}\label{sec1}

All graphs considered in this paper are finite, simple, and undirected.
For terms and symbols not defined in this paper, we refer the reader to \cite{D}.

Let $G$ be a graph.
Let $V(G)$ and $E(G)$ denote the {\it vertex set} and the {\it edge set} of $G$, respectively.
For a vertex $x\in V(G)$, let $N_{G}(x)$ denote the {\it neighborhood} of $x$ in $G$; thus $N_{G}(x)=\{y\in V(G): xy\in E(G)\}$.
For a subset $X$ of $V(G)$, let $N_{G}(X)=(\bigcup _{x\in X}N_{G}(x))\setminus X$, and let $G[X]$ (resp. $G-X$) denote the subgraph of $G$ induced by $X$ (resp. $V(G)\setminus X$).
Let $\alpha (G)$ denote the {\it independence number} of $G$, i.e., the maximum cardinality of an independent set of $G$.
Let $K_{n}$, $P_{n}$ and $K_{1,n}$ denote the {\it complete graph} of order $n$, the {\it path} of order $n$ and the {\it star} of order $n+1$, respectively.
For two positive integers $n_{1}$ and $n_{2}$, the {\it Ramsey number} $R(n_{1},n_{2})$ is the minimum positive integer $R$ such that any graph of order at least $R$ contains a clique of cardinality $n_{1}$ or an independent set of cardinality $n_{2}$.

For two graphs $G$ and $H$, $G$ is said to be {\it $H$-free} if $G$ contains no induced copy of $H$.
For a family $\HH$ of graphs, a graph $G$ is said to be {\it $\HH$-free} if $G$ is $H$-free for every $H\in \HH$.
In this context, the members of $\HH$ are called {\it forbidden subgraphs}.
For two families $\mathcal{H}_{1}$ and $\mathcal{H}_{2}$ of graphs, we write $\mathcal{H}_{1}\leq \mathcal{H}_{2}$ if for every $H_{2}\in \mathcal{H}_{2}$, there exists $H_{1}\in \mathcal{H}_{1}$ such that $H_{1}$ is an induced subgraph of $H_{2}$.
The relation ``$\leq $'' between two families of forbidden subgraphs was introduced in \cite{FKLOPS}.
Note that if $\mathcal{H}_{1}\leq \mathcal{H}_{2}$, then every $\mathcal{H}_{1}$-free graph is also $\mathcal{H}_{2}$-free.

Let $\AA$ be a family of graphs.
A family $\PP$ of subgraphs of $G$ is called an {\it $\AA$-cover} of $G$ if $\bigcup _{P\in \PP}V(P)=V(G)$ and each element of $\PP$ is isomorphic to a graph belonging to $\AA$.
Note that some elements of an $\AA$-cover of $G$ might have common vertices.
An $\AA$-cover $\PP$ of $G$ is called an {\it $\AA$-partition} of $G$ if the elements of $\PP$ are pairwise vertex-disjoint.
A $\{P_{i}:i\geq 1\}$-cover (resp. a $\{P_{i}:i\geq 1\}$-partition) of $G$ is called a {\it path cover} (resp. a {\it path partition}) of $G$.
Since $\{G[\{x\}]:x\in V(G)\}$ is a path partition of $G$ (and so a path cover of $G$), the minimum cardinality of a path cover (or a path partition) of any graph is well-defined.
The value $\min \{|\PP |:\PP\mbox{ is a path cover of }G\}$ (resp. $\min \{|\PP |:\PP\mbox{ is a path partition of }G\}$), denoted by ${\rm pc}(G)$ (resp. ${\rm pp}(G)$), is called the {\it path cover number} (resp. the {\it path partition number}) of $G$.
It is trivial that ${\rm pc}(G)\leq {\rm pp}(G)$.
Since a graph $G$ has a Hamiltonian path if and only if ${\rm pp}(G)=1$, the decision problem for the path partition number is a natural generalization of the Hamiltonian path problem.
In fact, it has been widely studied in, for example, \cite{H,I,MM,R,Y}.
Throughout this paper, we implicitly use the following fact.

\begin{fact}
\label{fact1}
Let $G$ be a graph, and let $\{X_{1},X_{2},\ldots ,X_{m}\}$ be a partition of $V(G)$.
Then ${\rm pc}(G)\leq \sum _{1\leq i\leq m}{\rm pc}(G[X_{i}])$ and ${\rm pp}(G)\leq \sum _{1\leq i\leq m}{\rm pp}(G[X_{i}])$.
\end{fact}

In this paper, we focus on the following conditions concerning a family $\HH$ of forbidden subgraphs:
\begin{enumerate}
\item[{\bf (A1)}]
There exists a constant $c_{1}=c_{1}(\HH)$ such that ${\rm pc}(G)\leq c_{1}$ for every connected $\HH$-free graph $G$.
\item[{\bf (A2)}]
There exists a constant $c_{2}=c_{2}(\HH)$ such that ${\rm pp}(G)\leq c_{2}$ for every connected $\HH$-free graph $G$.
\end{enumerate}
Our main aim is to characterize the finite families $\HH $ of connected graphs satisfying (A1) or (A2).

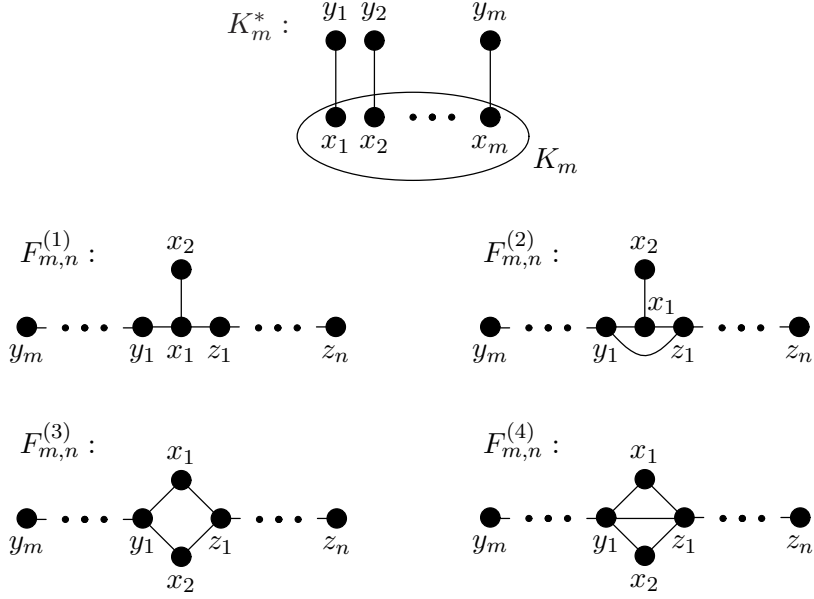
\begin{figure}
\begin{center}
{\unitlength 0.1in%
\begin{picture}(45.7500,29.8500)(-1.2000,-38.7500)%
\put(17.6000,-11.0000){\makebox(0,0)[rb]{$K^{*}_{m}:$}}%
%
\special{sh 1.000}%
\special{ia 2000 1105 50 50 0.0000000 6.2831853}%
\special{pn 8}%
\special{ar 2000 1105 50 50 0.0000000 6.2831853}%
%
\special{pn 8}%
\special{ar 2400 1605 600 230 0.0000000 6.2831853}%
%
\special{sh 1.000}%
\special{ia 2000 1505 50 50 0.0000000 6.2831853}%
\special{pn 8}%
\special{ar 2000 1505 50 50 0.0000000 6.2831853}%
%
\special{sh 1.000}%
\special{ia 2200 1505 50 50 0.0000000 6.2831853}%
\special{pn 8}%
\special{ar 2200 1505 50 50 0.0000000 6.2831853}%
%
\special{sh 1.000}%
\special{ia 2800 1505 50 50 0.0000000 6.2831853}%
\special{pn 8}%
\special{ar 2800 1505 50 50 0.0000000 6.2831853}%
%
\special{pn 4}%
\special{sh 1}%
\special{ar 2400 1505 16 16 0 6.2831853}%
\special{sh 1}%
\special{ar 2600 1505 16 16 0 6.2831853}%
\special{sh 1}%
\special{ar 2500 1505 16 16 0 6.2831853}%
\special{sh 1}%
\special{ar 2500 1505 16 16 0 6.2831853}%
%
\special{sh 1.000}%
\special{ia 2200 1105 50 50 0.0000000 6.2831853}%
\special{pn 8}%
\special{ar 2200 1105 50 50 0.0000000 6.2831853}%
%
\special{sh 1.000}%
\special{ia 2800 1105 50 50 0.0000000 6.2831853}%
\special{pn 8}%
\special{ar 2800 1105 50 50 0.0000000 6.2831853}%
%
\special{pn 8}%
\special{pa 2800 1105}%
\special{pa 2800 1505}%
\special{fp}%
\special{pa 2200 1505}%
\special{pa 2200 1105}%
\special{fp}%
\special{pa 2000 1105}%
\special{pa 2000 1505}%
\special{fp}%
\put(20.0000,-16.4500){\makebox(0,0){$x_{1}$}}%
\put(22.0000,-16.4500){\makebox(0,0){$x_{2}$}}%
\put(28.0000,-16.4500){\makebox(0,0){$x_{m}$}}%
\put(28.0000,-9.5500){\makebox(0,0){$y_{m}$}}%
\put(22.0000,-9.5500){\makebox(0,0){$y_{2}$}}%
\put(20.0000,-9.5500){\makebox(0,0){$y_{1}$}}%
\put(31.4000,-17.3000){\makebox(0,0){$K_{m}$}}%
%
\special{sh 1.000}%
\special{ia 400 2600 50 50 0.0000000 6.2831853}%
\special{pn 8}%
\special{ar 400 2600 50 50 0.0000000 6.2831853}%
%
\special{sh 1.000}%
\special{ia 1000 2600 50 50 0.0000000 6.2831853}%
\special{pn 8}%
\special{ar 1000 2600 50 50 0.0000000 6.2831853}%
%
\special{sh 1.000}%
\special{ia 1200 2600 50 50 0.0000000 6.2831853}%
\special{pn 8}%
\special{ar 1200 2600 50 50 0.0000000 6.2831853}%
%
\special{sh 1.000}%
\special{ia 1400 2600 50 50 0.0000000 6.2831853}%
\special{pn 8}%
\special{ar 1400 2600 50 50 0.0000000 6.2831853}%
%
\special{sh 1.000}%
\special{ia 2000 2600 50 50 0.0000000 6.2831853}%
\special{pn 8}%
\special{ar 2000 2600 50 50 0.0000000 6.2831853}%
%
\special{sh 1.000}%
\special{ia 1200 2300 50 50 0.0000000 6.2831853}%
\special{pn 8}%
\special{ar 1200 2300 50 50 0.0000000 6.2831853}%
%
\special{pn 4}%
\special{sh 1}%
\special{ar 600 2605 16 16 0 6.2831853}%
\special{sh 1}%
\special{ar 800 2605 16 16 0 6.2831853}%
\special{sh 1}%
\special{ar 700 2605 16 16 0 6.2831853}%
\special{sh 1}%
\special{ar 700 2605 16 16 0 6.2831853}%
%
\special{pn 4}%
\special{sh 1}%
\special{ar 1600 2605 16 16 0 6.2831853}%
\special{sh 1}%
\special{ar 1800 2605 16 16 0 6.2831853}%
\special{sh 1}%
\special{ar 1700 2605 16 16 0 6.2831853}%
\special{sh 1}%
\special{ar 1700 2605 16 16 0 6.2831853}%
%
\special{pn 8}%
\special{pa 400 2605}%
\special{pa 500 2605}%
\special{fp}%
\special{pa 900 2605}%
\special{pa 1000 2605}%
\special{fp}%
%
\special{pn 8}%
\special{pa 1400 2600}%
\special{pa 1500 2600}%
\special{fp}%
\special{pa 1900 2600}%
\special{pa 2000 2600}%
\special{fp}%
%
\special{pn 8}%
\special{pa 1400 2600}%
\special{pa 1000 2600}%
\special{fp}%
\special{pa 1200 2600}%
\special{pa 1200 2300}%
\special{fp}%
\put(4.0000,-27.4000){\makebox(0,0){$y_{m}$}}%
\put(10.0000,-27.4000){\makebox(0,0){$y_{1}$}}%
\put(12.0000,-27.4000){\makebox(0,0){$x_{1}$}}%
\put(14.0000,-27.4000){\makebox(0,0){$z_{1}$}}%
\put(20.0000,-27.4000){\makebox(0,0){$z_{n}$}}%
\put(12.0000,-21.7000){\makebox(0,0){$x_{2}$}}%
%
\special{sh 1.000}%
\special{ia 400 3600 50 50 0.0000000 6.2831853}%
\special{pn 8}%
\special{ar 400 3600 50 50 0.0000000 6.2831853}%
%
\special{sh 1.000}%
\special{ia 1000 3600 50 50 0.0000000 6.2831853}%
\special{pn 8}%
\special{ar 1000 3600 50 50 0.0000000 6.2831853}%
%
\special{sh 1.000}%
\special{ia 1405 3600 50 50 0.0000000 6.2831853}%
\special{pn 8}%
\special{ar 1405 3600 50 50 0.0000000 6.2831853}%
%
\special{sh 1.000}%
\special{ia 2005 3600 50 50 0.0000000 6.2831853}%
\special{pn 8}%
\special{ar 2005 3600 50 50 0.0000000 6.2831853}%
%
\special{pn 4}%
\special{sh 1}%
\special{ar 600 3605 16 16 0 6.2831853}%
\special{sh 1}%
\special{ar 800 3605 16 16 0 6.2831853}%
\special{sh 1}%
\special{ar 700 3605 16 16 0 6.2831853}%
\special{sh 1}%
\special{ar 700 3605 16 16 0 6.2831853}%
%
\special{pn 4}%
\special{sh 1}%
\special{ar 1605 3605 16 16 0 6.2831853}%
\special{sh 1}%
\special{ar 1805 3605 16 16 0 6.2831853}%
\special{sh 1}%
\special{ar 1705 3605 16 16 0 6.2831853}%
\special{sh 1}%
\special{ar 1705 3605 16 16 0 6.2831853}%
%
\special{pn 8}%
\special{pa 400 3605}%
\special{pa 500 3605}%
\special{fp}%
\special{pa 900 3605}%
\special{pa 1000 3605}%
\special{fp}%
%
\special{pn 8}%
\special{pa 1405 3600}%
\special{pa 1505 3600}%
\special{fp}%
\special{pa 1905 3600}%
\special{pa 2005 3600}%
\special{fp}%
\put(4.0000,-37.4000){\makebox(0,0){$y_{m}$}}%
\put(10.0000,-37.4000){\makebox(0,0){$y_{1}$}}%
\put(14.0500,-37.4000){\makebox(0,0){$z_{1}$}}%
\put(20.0500,-37.4000){\makebox(0,0){$z_{n}$}}%
%
\special{sh 1.000}%
\special{ia 1200 3400 50 50 0.0000000 6.2831853}%
\special{pn 8}%
\special{ar 1200 3400 50 50 0.0000000 6.2831853}%
%
\special{sh 1.000}%
\special{ia 1200 3800 50 50 0.0000000 6.2831853}%
\special{pn 8}%
\special{ar 1200 3800 50 50 0.0000000 6.2831853}%
%
\special{pn 8}%
\special{pa 1200 3800}%
\special{pa 1400 3600}%
\special{fp}%
\special{pa 1400 3600}%
\special{pa 1200 3400}%
\special{fp}%
\special{pa 1200 3400}%
\special{pa 1000 3600}%
\special{fp}%
\special{pa 1000 3600}%
\special{pa 1200 3800}%
\special{fp}%
\put(12.0000,-32.7000){\makebox(0,0){$x_{1}$}}%
\put(12.0000,-39.4000){\makebox(0,0){$x_{2}$}}%
%
\special{sh 1.000}%
\special{ia 2800 3595 50 50 0.0000000 6.2831853}%
\special{pn 8}%
\special{ar 2800 3595 50 50 0.0000000 6.2831853}%
%
\special{sh 1.000}%
\special{ia 3400 3595 50 50 0.0000000 6.2831853}%
\special{pn 8}%
\special{ar 3400 3595 50 50 0.0000000 6.2831853}%
%
\special{sh 1.000}%
\special{ia 3805 3595 50 50 0.0000000 6.2831853}%
\special{pn 8}%
\special{ar 3805 3595 50 50 0.0000000 6.2831853}%
%
\special{sh 1.000}%
\special{ia 4405 3595 50 50 0.0000000 6.2831853}%
\special{pn 8}%
\special{ar 4405 3595 50 50 0.0000000 6.2831853}%
%
\special{pn 4}%
\special{sh 1}%
\special{ar 3000 3600 16 16 0 6.2831853}%
\special{sh 1}%
\special{ar 3200 3600 16 16 0 6.2831853}%
\special{sh 1}%
\special{ar 3100 3600 16 16 0 6.2831853}%
\special{sh 1}%
\special{ar 3100 3600 16 16 0 6.2831853}%
%
\special{pn 4}%
\special{sh 1}%
\special{ar 4005 3600 16 16 0 6.2831853}%
\special{sh 1}%
\special{ar 4205 3600 16 16 0 6.2831853}%
\special{sh 1}%
\special{ar 4105 3600 16 16 0 6.2831853}%
\special{sh 1}%
\special{ar 4105 3600 16 16 0 6.2831853}%
%
\special{pn 8}%
\special{pa 2800 3600}%
\special{pa 2900 3600}%
\special{fp}%
\special{pa 3300 3600}%
\special{pa 3400 3600}%
\special{fp}%
%
\special{pn 8}%
\special{pa 3805 3595}%
\special{pa 3905 3595}%
\special{fp}%
\special{pa 4305 3595}%
\special{pa 4405 3595}%
\special{fp}%
\put(28.0000,-37.3500){\makebox(0,0){$y_{m}$}}%
\put(34.0000,-37.3500){\makebox(0,0){$y_{1}$}}%
\put(38.0500,-37.3500){\makebox(0,0){$z_{1}$}}%
\put(44.0500,-37.3500){\makebox(0,0){$z_{n}$}}%
%
\special{sh 1.000}%
\special{ia 3600 3395 50 50 0.0000000 6.2831853}%
\special{pn 8}%
\special{ar 3600 3395 50 50 0.0000000 6.2831853}%
%
\special{sh 1.000}%
\special{ia 3600 3795 50 50 0.0000000 6.2831853}%
\special{pn 8}%
\special{ar 3600 3795 50 50 0.0000000 6.2831853}%
%
\special{pn 8}%
\special{pa 3600 3795}%
\special{pa 3800 3595}%
\special{fp}%
\special{pa 3800 3595}%
\special{pa 3600 3395}%
\special{fp}%
\special{pa 3600 3395}%
\special{pa 3400 3595}%
\special{fp}%
\special{pa 3400 3595}%
\special{pa 3600 3795}%
\special{fp}%
\put(36.0000,-32.6500){\makebox(0,0){$x_{1}$}}%
\put(36.0000,-39.3500){\makebox(0,0){$x_{2}$}}%
%
\special{pn 8}%
\special{pa 3400 3600}%
\special{pa 3800 3600}%
\special{fp}%
\put(7.6000,-23.0000){\makebox(0,0)[rb]{$F^{(1)}_{m,n}:$}}%
\put(31.6000,-33.0000){\makebox(0,0)[rb]{$F^{(4)}_{m,n}:$}}%
\put(7.6000,-33.0000){\makebox(0,0)[rb]{$F^{(3)}_{m,n}:$}}%
%
\special{sh 1.000}%
\special{ia 2800 2600 50 50 0.0000000 6.2831853}%
\special{pn 8}%
\special{ar 2800 2600 50 50 0.0000000 6.2831853}%
%
\special{sh 1.000}%
\special{ia 3400 2600 50 50 0.0000000 6.2831853}%
\special{pn 8}%
\special{ar 3400 2600 50 50 0.0000000 6.2831853}%
%
\special{sh 1.000}%
\special{ia 3600 2600 50 50 0.0000000 6.2831853}%
\special{pn 8}%
\special{ar 3600 2600 50 50 0.0000000 6.2831853}%
%
\special{sh 1.000}%
\special{ia 3800 2600 50 50 0.0000000 6.2831853}%
\special{pn 8}%
\special{ar 3800 2600 50 50 0.0000000 6.2831853}%
%
\special{sh 1.000}%
\special{ia 4400 2600 50 50 0.0000000 6.2831853}%
\special{pn 8}%
\special{ar 4400 2600 50 50 0.0000000 6.2831853}%
%
\special{sh 1.000}%
\special{ia 3600 2300 50 50 0.0000000 6.2831853}%
\special{pn 8}%
\special{ar 3600 2300 50 50 0.0000000 6.2831853}%
%
\special{pn 4}%
\special{sh 1}%
\special{ar 3000 2605 16 16 0 6.2831853}%
\special{sh 1}%
\special{ar 3200 2605 16 16 0 6.2831853}%
\special{sh 1}%
\special{ar 3100 2605 16 16 0 6.2831853}%
\special{sh 1}%
\special{ar 3100 2605 16 16 0 6.2831853}%
%
\special{pn 4}%
\special{sh 1}%
\special{ar 4000 2605 16 16 0 6.2831853}%
\special{sh 1}%
\special{ar 4200 2605 16 16 0 6.2831853}%
\special{sh 1}%
\special{ar 4100 2605 16 16 0 6.2831853}%
\special{sh 1}%
\special{ar 4100 2605 16 16 0 6.2831853}%
%
\special{pn 8}%
\special{pa 2800 2605}%
\special{pa 2900 2605}%
\special{fp}%
\special{pa 3300 2605}%
\special{pa 3400 2605}%
\special{fp}%
%
\special{pn 8}%
\special{pa 3800 2600}%
\special{pa 3900 2600}%
\special{fp}%
\special{pa 4300 2600}%
\special{pa 4400 2600}%
\special{fp}%
%
\special{pn 8}%
\special{pa 3800 2600}%
\special{pa 3400 2600}%
\special{fp}%
\special{pa 3600 2600}%
\special{pa 3600 2300}%
\special{fp}%
\put(28.0000,-27.4000){\makebox(0,0){$y_{m}$}}%
\put(36.9000,-24.9000){\makebox(0,0){$x_{1}$}}%
\put(44.0000,-27.4000){\makebox(0,0){$z_{n}$}}%
\put(31.6000,-23.0000){\makebox(0,0)[rb]{$F^{(2)}_{m,n}:$}}%
\put(36.0000,-21.7000){\makebox(0,0){$x_{2}$}}%
%
\special{pn 8}%
\special{pa 3400 2600}%
\special{pa 3426 2629}%
\special{pa 3451 2656}%
\special{pa 3477 2682}%
\special{pa 3502 2705}%
\special{pa 3528 2724}%
\special{pa 3554 2739}%
\special{pa 3579 2748}%
\special{pa 3605 2750}%
\special{pa 3630 2745}%
\special{pa 3656 2734}%
\special{pa 3682 2718}%
\special{pa 3707 2697}%
\special{pa 3733 2673}%
\special{pa 3758 2646}%
\special{pa 3784 2618}%
\special{pa 3800 2600}%
\special{fp}%
\put(38.0000,-27.4000){\makebox(0,0){$z_{1}$}}%
\put(34.0000,-27.4000){\makebox(0,0){$y_{1}$}}%
\end{picture}}%

\caption{Graphs $K^{*}_{m}$, $F^{(1)}_{m,n}$, $F^{(2)}_{m,n}$, $F^{(3)}_{m,n}$ and $F^{(4)}_{m,n}$}
\label{f1}
\end{center}
\end{figure}

Let $m$ and $n$ be two positive integers.
We define five graphs which will be used as forbidden subgraphs in our main result (see Figure~\ref{f1}).
\begin{enumerate}[{$\bullet $}]
\item
Let $K^{*}_{m}$ denote the graph with $V(K^{*}_{m})=\{x_{i},y_{i}:1\leq i\leq m\}$ and $E(K^{*}_{m})=\{x_{i}x_{j}:1\leq i<j\leq m\}\cup \{x_{i}y_{i}:1\leq i\leq m\}$.
\item
Let $A:=\{x_{1},x_{2}\}\cup \{y_{i}:1\leq i\leq m\}\cup \{z_{i}:1\leq i\leq n\}$.
We define four graphs as follows:
\begin{enumerate}[{$\circ$}]
\item
Let $F^{(1)}_{m,n}$ denote the graph on $A$ such that $E(F^{(1)}_{m,n})=\{x_{1}x_{2},x_{1}y_{1},x_{1}z_{1}\}\cup \{y_{i}y_{i+1}:1\leq i\leq m-1\}\cup \{z_{i}z_{i+1}:1\leq i\leq n-1\}$.
\item
Let $F^{(2)}_{m,n}$ is the graph obtained from $F^{(1)}_{m,n}$ by adding the edge $y_{1}z_{1}$.
\item
Let $F^{(3)}_{m,n}$ denote the graph on $A$ such that $E(F^{(3)}_{m,n})=\{x_{1}y_{1},x_{1}z_{1},x_{2}y_{1},x_{2}z_{1}\}\cup \{y_{i}y_{i+1}:1\leq i\leq m-1\}\cup \{z_{i}z_{i+1}:1\leq i\leq n-1\}$.
\item
Let $F^{(4)}_{m,n}$ is the graph obtained from $F^{(3)}_{m,n}$ by adding the edge $y_{1}z_{1}$.
\end{enumerate}
\end{enumerate}
Our first main result is the following, which is proved in Section~\ref{sec2}.

\begin{thm}
\label{mainthm}
Let $\HH$ be a finite family of connected graphs.
Then the following hold:
\begin{enumerate}
\item[{\upshape(i)}]
The family $\HH$ satisfies (A1) if and only if $\HH\leq \{K_{1,n},K^{*}_{n},F^{(1)}_{n,n},F^{(2)}_{n,n}\}$ for an integer $n\geq 2$.
\item[{\upshape(ii)}]
The family $\HH$ satisfies (A2) if and only if $\HH\leq \{K_{1,n},K^{*}_{n},F^{(1)}_{n,n},F^{(2)}_{n,n},F^{(3)}_{n,n},F^{(4)}_{n,n}\}$ for an integer $n\geq 2$.
\end{enumerate}
\end{thm}

Our motivation derives from two different lines of research.
The first one is forbidden subgraph conditions for the existence of a Hamiltonian path.
Now we focus on the condition that
\begin{align}
\mbox{every connected $\HH$-free graph (of sufficiently large order) has a Hamiltonian path}\label{cond-intro-1}
\end{align}
for a family $\HH$ of connected graphs.
Duffus et al.~\cite{DGJ} proved $\HH=\{K_{1,3},K^{*}_{3}\}$ satisfies (\ref{cond-intro-1}), and Faudree and Gould~\cite{FG} showed that if a family $\HH$ satisfying (\ref{cond-intro-1}) consists of two connected graphs, then $\HH\leq \{K_{1,3},K^{*}_{3}\}$.
Thereafter a series by Gould and Harris~\cite{GH1,GH2,GH3} characterized the families $\HH$ of connected graphs with $|\HH|=3$ satisfying (\ref{cond-intro-1}).
Since a graph has a Hamiltonian path if and only if its path cover number (or its path partition number) is exactly one, it is natural to study the forbidden subgraph conditions assuring us that the path cover/partition number is bounded by a constant as a next step.
Our main result gives a complete solution for the problem in a sense.

Our second motivation is an analysis of gap between minimum $\AA$-covers and minimum $\AA$-partitions.
A path cover/partition, which are main topic in this paper, is just one of examples of $\AA$-cover/partition problems, and there also exist many other cover/partition problems.
One of representative other examples is the case where $\AA$ is the family of all stars, where we regard $K_{1}$ as one of stars.
If we define the star cover number and the star partition number in the same way as ${\rm pc}(G)$ and ${\rm pp}(G)$, we can easily verify that the values are always equivalent.
(Indeed, the star cover number also equals to the {\it domination number}, which is one of classical invariants in graph theory.
The forbidden subgraph conditions assuring us that the domination number is bounded by a constant were characterized in \cite{F}.)
On the other hand, as it is evident from Theorem~\ref{mainthm}, there is a gap between the path cover number and the path partition number.
By Theorem~\ref{mainthm}, we discover that $F^{(3)}_{n,n}$ and $F^{(4)}_{n,n}$ play an important role for essential structures giving such a gap.

We also obtain an analogy of Theorem~\ref{mainthm} considering a cycle cover/partition problem.
A $\{K_{1},K_{2},C_{i}:i\geq 3\}$-cover (resp. a $\{K_{1},K_{2},C_{i}:i\geq 3\}$-partition) of $G$ is called a {\it cycle cover} (resp. a {\it cycle partition}) of $G$.
The value $\min \{|\PP |:\PP\mbox{ is a cycle cover of }G\}$ (resp. $\min \{|\PP |:\PP\mbox{ is a cycle partition of }G\}$), denoted by ${\rm cc}(G)$ (resp. ${\rm cp}(G)$), is called the {\it cycle cover number} (resp. the {\it cycle partition number}) of $G$.
Since trees (or graphs having a vertex of degree one) has no $\{C_{i}:i\geq 3\}$-cover, one sometimes focuses on cycle covers/partitions of general graphs instead of $\{C_{i}:i\geq 3\}$-covers/partitions (see, for example, \cite{EL,Fj}).
In Section~\ref{sec3}, as the second result, we characterize the families $\HH$ of forbidden subgraphs satisfying one of the following:
\begin{enumerate}
\item[{\bf (A'1)}]
There exists a constant $c_{1}=c_{1}(\HH)$ such that ${\rm cc}(G)\leq c_{1}$ for every connected $\HH$-free graph $G$.
\item[{\bf (A'2)}]
There exists a constant $c_{2}=c_{2}(\HH)$ such that ${\rm cp}(G)\leq c_{2}$ for every connected $\HH$-free graph $G$.
\end{enumerate}

\begin{thm}
\label{mainthm2}
Let $\HH$ be a family of connected graphs.
Then the following are equivalent.
\begin{enumerate}
\item[{\upshape(i)}]
The family $\HH$ satisfies (A'1).
\item[{\upshape(ii)}]
The family $\HH$ satisfies (A'2).
\item[{\upshape(iii)}]
For an integer $n\geq 2$, $\HH\leq \{K_{1,n},K^{*}_{n},P_{n}\}$.
\end{enumerate}
\end{thm}

We conclude this section by defining a new Ramsey-type concept concerning the path cover/partition number.
Let $\HH$ be a family of graphs.
The {\it path cover Ramsey number} $R^{\rm pc}(\HH)$ (resp. the {\it path partition Ramsey number} $R^{\rm pp}(\HH)$) is the minimum positive integer $R$ such that any connected graph $G$ with ${\rm pc}(G)\geq R$ (resp. ${\rm pp}(G)\geq R$) contains an induced copy of an element of $\HH$, where $R^{\rm pc}(\HH)=\infty $ (resp. $R^{\rm pp}(\HH)=\infty $) if such an integer does not exist.
Then it follows from Theorem~\ref{mainthm} that the following hold:
\begin{enumerate}
\item[{\bf (P1)}]
For a finite family $\HH$ of connected graphs, $R^{\rm pc}(\HH)$ is a finite number if and only if $\HH\leq \{K_{1,n},K^{*}_{n},F^{(1)}_{n,n},F^{(2)}_{n,n}\}$ for an integer $n\geq 2$.
\item[{\bf (P2)}]
For a finite family $\HH$ of connected graphs, $R^{\rm pp}(\HH)$ is a finite number if and only if $\HH\leq \{K_{1,n},K^{*}_{n},F^{(1)}_{n,n},F^{(2)}_{n,n},F^{(3)}_{n,n},F^{(4)}_{n,n}\}$ for an integer $n\geq 2$.
\end{enumerate}

Note that $R^{\rm pc}(\HH)=2$ if and only if $R^{\rm pp}(\HH)=2$.
As we mentioned above, it is known that $R^{\rm pc}(\{K_{1,3},K^{*}_{3}\})=2$ and the study of triples $\{H_{1},H_{2},H_{3}\}$ of connected graphs with $R^{\rm pc}(\{H_{1},H_{2},H_{3}\})=2$ is completed.
Since the $K_{1,3}$-freeness tends to give an important structure to many Hamiltonian properties, one might be interested in a relationship between such new Ramsey-type values and $K_{1,3}$-freeness.
Here we focus on the values $R^{\rm pc}(\HH)$ and $R^{\rm pp}(\HH)$ for the case where $\HH$ contains $K_{1,3}$.
Note that for positive integers $m$ and $n$ with $m+n\geq 3$, all of $F^{(1)}_{m,n}$, $F^{(3)}_{m,n}$ and $F^{(4)}_{m,n}$ contain $K_{1,3}$ as an induced copy.
Thus if $K_{1,3}\in \HH$, then
$$
R^{\rm pc}(\HH)=R^{\rm pc}(\HH\setminus \{F^{(1)}_{m,n},F^{(3)}_{m,n},F^{(4)}_{m,n}:m\geq 1,~n\geq 1,~m+n\geq 3\})
$$
and
$$
R^{\rm pp}(\HH)=R^{\rm pp}(\HH\setminus \{F^{(1)}_{m,n},F^{(3)}_{m,n},F^{(4)}_{m,n}:m\geq 1,~n\geq 1,~m+n\geq 3\}).
$$
Considering (P1) and (P2), we leave the following open problem which will be a next interesting target on this concept for readers.

\begin{problem}
\label{prob1}
For positive integers $p$, $q$ and $r$ with $p\geq 3$ and $q+r\geq 4$ and for a family $\HH$ of graphs with $\HH\leq \{K_{1,3},K^{*}_{p},F^{(2)}_{q,r}\}$, determine the value $R^{\rm pc}(\HH)$ and $R^{\rm pp}(\HH)$.
\end{problem}

\section{Proof of Theorem~\ref{mainthm}}\label{sec2}

\subsection{The ``if'' parts of Theorem~\ref{mainthm}}\label{sec2.2}

In this subsection, we prove the following theorem, which implies that the ``if'' parts of Theorem~\ref{mainthm} hold.

\begin{thm}
\label{thm1}
Let $n\geq 2$ be an integer.
Then the following hold:
\begin{enumerate}
\item[{\upshape(i)}]
There exists a constant $c_{1}=c_{1}(n)$ depending on $n$ only such that ${\rm pc}(G)\leq c_{1}$ for every connected $\{K_{1,n},K^{*}_{n},F^{(1)}_{n,n},F^{(2)}_{n,n}\}$-free graph $G$.
\item[{\upshape(ii)}]
There exists a constant $c_{2}=c_{2}(n)$ depending on $n$ only such that ${\rm pp}(G)\leq c_{2}$ for every connected $\{K_{1,n},K^{*}_{n},F^{(1)}_{n,n},F^{(2)}_{n,n},F^{(3)}_{n,n},F^{(4)}_{n,n}\}$-free graph $G$.
\end{enumerate}
\end{thm}

The following lemma is well-known (or it is also obtained from a result on digraph by Gallai and Milgram~\cite{GM}).
So many readers can skip are advised to skip the proof.

\begin{lem}
\label{lem2.0}
For a graph $G$, ${\rm pc}(G)\leq {\rm pp}(G)\leq \alpha (G)$.
\end{lem}
\proof
Since a path partition of $G$ is also a path cover of $G$, we have ${\rm pc}(G)\leq {\rm pp}(G)$.

Let $\PP$ be a path partition of $G$ with $|\PP|={\rm pp}(G)$, and write $\PP=\{Q_{i}:1\leq i\leq {\rm pp}(G)\}$.
For each $i$ with $1\leq i\leq {\rm pp}(G)$, let $x_{i}$ be an endvertex of $Q_{i}$.
If $x_{i}x_{j}\in E(G)$ for some $1\leq i<j\leq {\rm pp}(G)$, then the graph $Q$ obtained from $Q_{i}$ and $Q_{j}$ by joining the edge $x_{i}x_{j}$ is a path, and hence $\PP'=(\PP\setminus \{Q_{i},Q_{j}\})\cup \{Q\}$ is a path partition of $G$ with $|\PP'|={\rm pp}(G)-1$, which contradicts the definition of the path partition number.
Thus $\{x_{i}:1\leq i\leq {\rm pp}(G)\}$ is an independent set of $G$, and hence ${\rm pp}(G)\leq \alpha (G)$.
\qed

\begin{lem}
\label{lem2.1}
Let $n\geq 2$ and $\alpha \geq 1$ be integers.
Let $G$ be a $\{K_{1,n},K^{*}_{n}\}$-free graph, and let $X$ be a subset of $V(G)$ with $\alpha (G[X])\leq \alpha $.
Then $\alpha (G[N_{G}(X)])\leq (n-1)R(n,\alpha +1)-1$.
\end{lem}
\proof
By way of contradiction, we suppose that there exists a subset $Y$ of $N_{G}(X)$ such that $Y$ is an independent set of $G$ and $|Y|=(n-1)R(n,\alpha +1)$.
Take a subset $X_{0}$ of $X$ with $Y\subseteq N_{G}(X_{0})$ so that $|X_{0}|$ is as small as possible.
If $|X_{0}|\leq R(n,\alpha +1)-1$, then $\frac{|Y|}{|X_{0}|}\geq \frac{(n-1)R(n,\alpha +1)}{R(n,\alpha +1)-1}>n-1$, and hence there exists a vertex $x_{0}\in X_{0}$ with $|N_{G}(x_{0})\cap Y|\geq n$, which contradicts the $K_{1,n}$-freeness of $G$.
Thus $|X_{0}|\geq R(n,\alpha +1)$.
Since $\alpha (G[X_{0}])\leq \alpha (G[X])\leq \alpha $, this implies that there exists a subset $X_{1}$ of $X_{0}$ such that $X_{1}$ is a clique of $G$ and $|X_{1}|=n$.
By the minimality of $X_{0}$, $(N_{G}(x)\cap Y)\setminus N_{G}(X_{0}\setminus \{x\})\neq \emptyset $ for every $x\in X_{0}$.
For each $x\in X_{0}$, let $y_{x}\in (N_{G}(x)\cap Y)\setminus N_{G}(X_{0}\setminus \{x\})$.
Then $X_{1}\cup \{y_{x}:x\in X_{1}\}$ induces a copy of $K^{*}_{n}$ in $G$, which contradicts the $K^{*}_{n}$-freeness of $G$.
\qed

\begin{figure}
\begin{center}
{\unitlength 0.1in%
\begin{picture}(31.9000,26.1000)(2.1000,-29.0000)%
%
\special{sh 1.000}%
\special{ia 1200 2200 50 50 0.0000000 6.2831853}%
\special{pn 8}%
\special{ar 1200 2200 50 50 0.0000000 6.2831853}%
%
\special{sh 1.000}%
\special{ia 1400 2200 50 50 0.0000000 6.2831853}%
\special{pn 8}%
\special{ar 1400 2200 50 50 0.0000000 6.2831853}%
%
\special{sh 1.000}%
\special{ia 2000 2200 50 50 0.0000000 6.2831853}%
\special{pn 8}%
\special{ar 2000 2200 50 50 0.0000000 6.2831853}%
%
\special{sh 1.000}%
\special{ia 2600 2200 50 50 0.0000000 6.2831853}%
\special{pn 8}%
\special{ar 2600 2200 50 50 0.0000000 6.2831853}%
%
\special{sh 1.000}%
\special{ia 2600 2800 50 50 0.0000000 6.2831853}%
\special{pn 8}%
\special{ar 2600 2800 50 50 0.0000000 6.2831853}%
%
\special{sh 1.000}%
\special{ia 2000 2800 50 50 0.0000000 6.2831853}%
\special{pn 8}%
\special{ar 2000 2800 50 50 0.0000000 6.2831853}%
%
\special{sh 1.000}%
\special{ia 1400 2800 50 50 0.0000000 6.2831853}%
\special{pn 8}%
\special{ar 1400 2800 50 50 0.0000000 6.2831853}%
%
\special{sh 1.000}%
\special{ia 1200 2800 50 50 0.0000000 6.2831853}%
\special{pn 8}%
\special{ar 1200 2800 50 50 0.0000000 6.2831853}%
%
\special{pn 8}%
\special{pa 1100 2100}%
\special{pa 2200 2100}%
\special{pa 2200 2900}%
\special{pa 1100 2900}%
\special{pa 1100 2100}%
\special{pa 2200 2100}%
\special{fp}%
%
\special{pn 8}%
\special{pa 1200 2200}%
\special{pa 1500 2200}%
\special{fp}%
\special{pa 1900 2200}%
\special{pa 2000 2200}%
\special{fp}%
%
\special{pn 8}%
\special{pa 1200 2800}%
\special{pa 1500 2800}%
\special{fp}%
\special{pa 1900 2800}%
\special{pa 2000 2800}%
\special{fp}%
%
\special{pn 4}%
\special{sh 1}%
\special{ar 1600 2815 16 16 0 6.2831853}%
\special{sh 1}%
\special{ar 1800 2815 16 16 0 6.2831853}%
\special{sh 1}%
\special{ar 1700 2815 16 16 0 6.2831853}%
\special{sh 1}%
\special{ar 1700 2815 16 16 0 6.2831853}%
%
\special{pn 4}%
\special{sh 1}%
\special{ar 1600 2215 16 16 0 6.2831853}%
\special{sh 1}%
\special{ar 1800 2215 16 16 0 6.2831853}%
\special{sh 1}%
\special{ar 1700 2215 16 16 0 6.2831853}%
\special{sh 1}%
\special{ar 1700 2215 16 16 0 6.2831853}%
%
\special{pn 4}%
\special{sh 1}%
\special{ar 2800 2500 16 16 0 6.2831853}%
\special{sh 1}%
\special{ar 2800 2400 16 16 0 6.2831853}%
\special{sh 1}%
\special{ar 2800 2600 16 16 0 6.2831853}%
\special{sh 1}%
\special{ar 2800 2600 16 16 0 6.2831853}%
%
\special{pn 8}%
\special{pa 2600 2200}%
\special{pa 2634 2205}%
\special{pa 2668 2211}%
\special{pa 2699 2220}%
\special{pa 2727 2233}%
\special{pa 2751 2251}%
\special{pa 2771 2275}%
\special{pa 2787 2303}%
\special{pa 2800 2330}%
\special{fp}%
%
\special{pn 8}%
\special{pa 2600 2800}%
\special{pa 2634 2795}%
\special{pa 2668 2789}%
\special{pa 2699 2780}%
\special{pa 2727 2767}%
\special{pa 2751 2749}%
\special{pa 2771 2725}%
\special{pa 2787 2697}%
\special{pa 2800 2670}%
\special{fp}%
\put(16.0000,-25.0000){\makebox(0,0){$X_{0}$}}%
%
\special{pn 8}%
\special{pa 2500 2100}%
\special{pa 2900 2100}%
\special{pa 2900 2900}%
\special{pa 2500 2900}%
\special{pa 2500 2100}%
\special{pa 2900 2100}%
\special{fp}%
%
\special{pn 8}%
\special{pa 2000 2200}%
\special{pa 2600 2200}%
\special{fp}%
\special{pa 2600 2800}%
\special{pa 2000 2800}%
\special{fp}%
\put(23.5000,-29.2000){\makebox(0,0){$P$}}%
\put(12.1000,-23.3000){\makebox(0,0){$u_{1}$}}%
\put(20.0000,-23.3000){\makebox(0,0){$u_{n_{0}}$}}%
\put(19.2000,-26.7000){\makebox(0,0){$u_{m-n_{0}+1}$}}%
\put(12.2000,-26.7000){\makebox(0,0){$u_{m}$}}%
%
\special{pn 8}%
\special{ar 3000 1670 400 200 0.0000000 6.2831853}%
%
\special{pn 20}%
\special{pa 2750 2100}%
\special{pa 3000 1870}%
\special{fp}%
%
\special{pn 8}%
\special{ar 1500 1230 400 150 0.0000000 6.2831853}%
%
\special{pn 8}%
\special{ar 1500 830 400 150 0.0000000 6.2831853}%
\put(15.0000,-12.3000){\makebox(0,0){$X_{2}$}}%
\put(15.0000,-8.3000){\makebox(0,0){$X_{3}$}}%
%
\special{pn 20}%
\special{pa 1500 1080}%
\special{pa 1500 980}%
\special{fp}%
%
\special{pn 20}%
\special{pa 1500 680}%
\special{pa 1500 580}%
\special{fp}%
%
\special{pn 4}%
\special{sh 1}%
\special{ar 1500 390 16 16 0 6.2831853}%
\special{sh 1}%
\special{ar 1500 290 16 16 0 6.2831853}%
\special{sh 1}%
\special{ar 1500 490 16 16 0 6.2831853}%
\special{sh 1}%
\special{ar 1500 490 16 16 0 6.2831853}%
%
\special{pn 20}%
\special{pa 1500 1480}%
\special{pa 1500 1380}%
\special{fp}%
\put(4.1000,-16.8000){\makebox(0,0){$X_{1}$}}%
\put(30.0000,-16.7000){\makebox(0,0){$Y$}}%
%
\special{pn 8}%
\special{ar 1000 1680 400 150 0.0000000 6.2831853}%
%
\special{pn 8}%
\special{ar 2000 1680 400 150 0.0000000 6.2831853}%
%
\special{pn 20}%
\special{pa 2660 2100}%
\special{pa 2060 1830}%
\special{fp}%
\special{pa 1950 1830}%
\special{pa 1650 2100}%
\special{fp}%
\special{pa 1550 2100}%
\special{pa 1000 1830}%
\special{fp}%
%
\special{pn 8}%
\special{pa 550 1480}%
\special{pa 2450 1480}%
\special{pa 2450 1880}%
\special{pa 550 1880}%
\special{pa 550 1480}%
\special{pa 2450 1480}%
\special{fp}%
\end{picture}}%

\caption{Path $P$ and sets $X_{i}$ and $Y$}
\label{f2}
\end{center}
\end{figure}
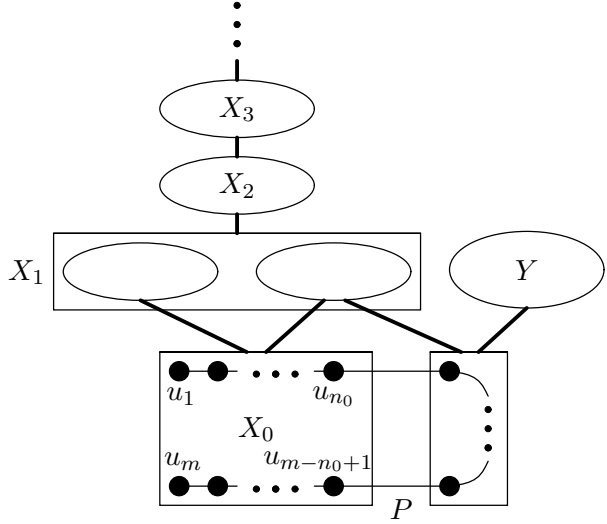

In the remainder of this subsection, we fix an integer $n\geq 2$ and a connected $\{K_{1,n},K^{*}_{n},F^{(1)}_{n,n},F^{(2)}_{n,n}\}$-free graph $G$.
Set $n_{0}=\max\{\lceil \frac{n^{2}-n-2}{2}\rceil ,n\}$.
Take a longest induced path $P$ of $G$, and write $P=u_{1}u_{2}\cdots u_{m}$.
Let $X_{0}=\{u_{i}:1\leq i\leq n_{0}\mbox{ or }m-n_{0}+1\leq i\leq m\}$ and $Y=N_{G}(V(P)\setminus X_{0})\setminus (X_{0}\cup N_{G}(X_{0}))$.
Note that if $|V(P)|\leq 2n_{0}$, then $X_{0}=V(P)$ and $Y=\emptyset $.
We further remark that $N_{G}(y)\cap V(P)\subseteq \{u_{i}:n_{0}+1\leq i\leq m-n_{0}\}$ for every $y\in Y$ (and in the remainder of this subsection, we frequently use the fact without mentioning).
For each $i$ with $n_{0}+1\leq i\leq m-n_{0}$, let $Y_{i}=\{y\in Y:\min\{j:n_{0}+1\leq j\leq m-n_{0},~yu_{j}\in E(G)\}=i\}$.
Now we recursively define the sets $X_{i}~(i\geq 1)$ as follows:
Let $X_{1}=N_{G}(X_{0})\setminus V(P)$, and for $i$ with $i\geq 2$, let $X_{i}=N_{G}(X_{i-1})\setminus (V(P)\cup Y\cup (\bigcup _{1\leq j\leq i-1}X_{j}))$ (see Figure~\ref{f2}).
Then $X_{1}\cap Y=\emptyset $ and $X_{1}\cup Y=N_{G}(V(P))$.

\begin{lem}
\label{lem2.2}
We have $X_{2n_{0}}=\emptyset $.
\end{lem}
\proof
Suppose that $X_{2n_{0}}\neq \emptyset $.
Let $x_{2n_{0}}\in X_{2n_{0}}$.
Then we can recursively take a vertex $x_{2n_{0}-i}\in N_{G}(x_{2n_{0}-i+1})\cap X_{2n_{0}-i}$ for $i$ with $1\leq i\leq 2n_{0}$.
Note that $x_{0}=u_{k}$ for some $k$ with $1\leq k\leq n_{0}$ or $m-n_{0}+1\leq k\leq m$.
By symmetry, we may assume that $1\leq k\leq n_{0}$.
Under this condition, we choose $k$ so that $k$ is as large as possible.
Since $x_{0}x_{1}\cdots x_{2n_{0}}$ is an induced path of $G$ having $2n_{0}+1$ vertices, it follows from the maximality of $P$ that $|V(P)|\geq 2n_{0}+1$.
In particular, $V(P)\setminus X_{0}\neq \emptyset $.

If $N_{G}(x_{1})\cap (V(P)\setminus X_{0})=\emptyset $, then $x_{2n_{0}}x_{2n_{0}-1}\cdots x_{1}u_{k}u_{k+1}\cdots u_{m-n_{0}}$ is an induced path of $G$ having $2n_{0}+m-n_{0}-k+1~(\geq m+1)$ vertices, which contradicts the maximality of $P$.
Thus $N_{G}(x_{1})\cap (V(P)\setminus X_{0})\neq \emptyset $.

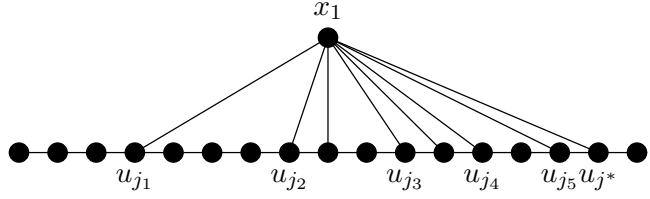
\begin{figure}
\begin{center}
{\unitlength 0.1in%
\begin{picture}(33.0000,8.8000)(7.5000,-22.7500)%
%
\special{sh 1.000}%
\special{ia 1000 2200 50 50 0.0000000 6.2831853}%
\special{pn 8}%
\special{ar 1000 2200 50 50 0.0000000 6.2831853}%
%
\special{sh 1.000}%
\special{ia 1200 2200 50 50 0.0000000 6.2831853}%
\special{pn 8}%
\special{ar 1200 2200 50 50 0.0000000 6.2831853}%
%
\special{sh 1.000}%
\special{ia 1400 2200 50 50 0.0000000 6.2831853}%
\special{pn 8}%
\special{ar 1400 2200 50 50 0.0000000 6.2831853}%
%
\special{sh 1.000}%
\special{ia 1600 2200 50 50 0.0000000 6.2831853}%
\special{pn 8}%
\special{ar 1600 2200 50 50 0.0000000 6.2831853}%
%
\special{sh 1.000}%
\special{ia 1800 2200 50 50 0.0000000 6.2831853}%
\special{pn 8}%
\special{ar 1800 2200 50 50 0.0000000 6.2831853}%
%
\special{sh 1.000}%
\special{ia 2000 2200 50 50 0.0000000 6.2831853}%
\special{pn 8}%
\special{ar 2000 2200 50 50 0.0000000 6.2831853}%
%
\special{sh 1.000}%
\special{ia 2200 2200 50 50 0.0000000 6.2831853}%
\special{pn 8}%
\special{ar 2200 2200 50 50 0.0000000 6.2831853}%
%
\special{sh 1.000}%
\special{ia 2400 2200 50 50 0.0000000 6.2831853}%
\special{pn 8}%
\special{ar 2400 2200 50 50 0.0000000 6.2831853}%
\put(14.0000,-23.4000){\makebox(0,0){$u_{j_{1}}$}}%
%
\special{sh 1.000}%
\special{ia 2600 2200 50 50 0.0000000 6.2831853}%
\special{pn 8}%
\special{ar 2600 2200 50 50 0.0000000 6.2831853}%
%
\special{sh 1.000}%
\special{ia 2800 2200 50 50 0.0000000 6.2831853}%
\special{pn 8}%
\special{ar 2800 2200 50 50 0.0000000 6.2831853}%
%
\special{sh 1.000}%
\special{ia 3000 2200 50 50 0.0000000 6.2831853}%
\special{pn 8}%
\special{ar 3000 2200 50 50 0.0000000 6.2831853}%
%
\special{sh 1.000}%
\special{ia 3200 2200 50 50 0.0000000 6.2831853}%
\special{pn 8}%
\special{ar 3200 2200 50 50 0.0000000 6.2831853}%
%
\special{sh 1.000}%
\special{ia 3400 2200 50 50 0.0000000 6.2831853}%
\special{pn 8}%
\special{ar 3400 2200 50 50 0.0000000 6.2831853}%
%
\special{sh 1.000}%
\special{ia 3600 2200 50 50 0.0000000 6.2831853}%
\special{pn 8}%
\special{ar 3600 2200 50 50 0.0000000 6.2831853}%
%
\special{sh 1.000}%
\special{ia 3800 2200 50 50 0.0000000 6.2831853}%
\special{pn 8}%
\special{ar 3800 2200 50 50 0.0000000 6.2831853}%
%
\special{sh 1.000}%
\special{ia 4000 2200 50 50 0.0000000 6.2831853}%
\special{pn 8}%
\special{ar 4000 2200 50 50 0.0000000 6.2831853}%
\put(22.0000,-23.4000){\makebox(0,0){$u_{j_{2}}$}}%
\put(28.0000,-23.4000){\makebox(0,0){$u_{j_{3}}$}}%
\put(32.0000,-23.4000){\makebox(0,0){$u_{j_{4}}$}}%
\put(36.0000,-23.4000){\makebox(0,0){$u_{j_{5}}$}}%
\put(38.0000,-23.4000){\makebox(0,0){$u_{j^{*}}$}}%
%
\special{sh 1.000}%
\special{ia 800 2200 50 50 0.0000000 6.2831853}%
\special{pn 8}%
\special{ar 800 2200 50 50 0.0000000 6.2831853}%
%
\special{sh 1.000}%
\special{ia 2400 1600 50 50 0.0000000 6.2831853}%
\special{pn 8}%
\special{ar 2400 1600 50 50 0.0000000 6.2831853}%
%
\special{pn 8}%
\special{pa 2400 1600}%
\special{pa 1400 2200}%
\special{fp}%
\special{pa 2200 2200}%
\special{pa 2400 1600}%
\special{fp}%
\special{pa 2400 1600}%
\special{pa 2400 2200}%
\special{fp}%
\special{pa 2400 1600}%
\special{pa 2800 2200}%
\special{fp}%
\special{pa 3000 2200}%
\special{pa 2400 1600}%
\special{fp}%
\special{pa 2400 1600}%
\special{pa 3200 2200}%
\special{fp}%
\special{pa 3600 2200}%
\special{pa 2400 1600}%
\special{fp}%
\special{pa 2400 1600}%
\special{pa 3800 2200}%
\special{fp}%
%
\special{pn 8}%
\special{pa 800 2200}%
\special{pa 4000 2200}%
\special{fp}%
\put(24.0000,-14.6000){\makebox(0,0){$x_{1}$}}%
\end{picture}}%

\caption{An example for $s=5$}
\label{f3}
\end{center}
\end{figure}

Now we consider an operation recursively defining integers $j_{1},j_{2},\ldots $ with $1\leq j_{p}\leq m~(p\geq 1)$ and $j_{1}<j_{2}<\cdots $ as follows (see Figure~\ref{f3}):
Let $j_{1}=\min\{j:1\leq j\leq m,~x_{1}u_{j}\in E(G)\}$.
For $p\geq 2$, we assume that the integer $j_{p-1}$ has defined.
If $\{j:j_{p-1}+2\leq j\leq m,~x_{1}u_{j}\in E(G)\}\neq \emptyset $, we let $j_{p}=\min\{j:j_{p-1}+2\leq j\leq m,~x_{1}u_{j}\in E(G)\}$; otherwise, we finish the operation.
Let $S=\{u_{j_{p}}:p\geq 1\}$, and set $s=|S|$.
Let $j^{*}=\max\{j:1\leq j\leq m,~x_{1}u_{j}\in E(G)\}$.
Note that $j^{*}\in \{j_{s},j_{s}+1\}$.
Since $j_{p}\geq j_{p-1}+2$, $S$ is an independent set of $G$.
Since $G$ is $K_{1,n}$-free and $\{x_{1},x_{2}\}\cup S$ induces a copy of $K_{1,s+1}$ in $G$, we have $s+1\leq n-1$.

For the moment, suppose that $s=1$.
Since $N_{G}(x_{1})\cap \{u_{j}:1\leq j\leq n_{0}\}\neq \emptyset $ and $N_{G}(x_{1})\cap (V(P)\setminus X_{0})\neq \emptyset $, this forces $N_{G}(x_{1})=\{u_{n_{0}},u_{n_{0}+1}\}$.
Then
$$
\{x_{1},x_{2},u_{n_{0}},u_{n_{0}-1},\ldots ,u_{n_{0}-n+1},u_{n_{0}+1},u_{n_{0}+2},\ldots ,u_{n_{0}+n}\}
$$
induces a copy of $F^{(2)}_{n,n}$ in $G$, which is a contradiction.
Thus $s\geq 2$.

Let $Q_{1}=u_{1}u_{2}\cdots u_{j_{1}}$ and $Q_{s+1}=u_{j^{*}}u_{j^{*}+1}\cdots u_{m}$ be subpaths of $P$.
For $p$ with $2\leq p\leq s$, let $Q_{p}=u_{j_{p-1}+2}u_{j_{p-1}+3}\cdots u_{j_{p}}$ be a subpath of $P$.
Then $V(P)\setminus (\bigcup _{1\leq p\leq s+1}V(Q_{p}))=\{u_{j_{p}+1}:1\leq p\leq s-1\}$, and hence
\begin{align*}
2n_{0}+1&\leq |V(P)|\\
&= \left|V(P)\setminus \left(\bigcup _{1\leq p\leq s+1}V(Q_{p}) \right) \right|+ \left|\bigcup _{1\leq p\leq s+1}V(Q_{p}) \right|\\
&\leq (s-1)+\sum _{1\leq p\leq s+1}|V(Q_{p})|\\
&= \sum _{1\leq p\leq s+1}(|V(Q_{p})|+1)-2.
\end{align*}
This implies that $\sum _{1\leq p\leq s+1}(|V(Q_{p})|+1)\geq 2n_{0}+3\geq 2\lceil \frac{n^{2}-n-2}{2}\rceil +3\geq n^{2}-n+1$.
If $|V(Q_{p})|\leq n-1$ for all $p$ with $1\leq p\leq s+1$, then $n^{2}-n+1\leq \sum _{1\leq p\leq s+1}(|V(Q_{p})|+1)\leq (s+1)n\leq (n-1)n$, which is a contradiction.
Thus $|V(Q_{q})|\geq n$ for some $q$ with $1\leq q\leq s+1$.

Note that $|N_{G}(x_{1})\cap V(Q_{q})|=1$.
Write $N_{G}(x_{1})\cap V(Q_{q})=\{u_{j}\}$.
If $q\neq s+1$, then $j\in \{j_{p}:1\leq p\leq s\}$; otherwise, $j=j^{*}~(\in \{j_{s},j_{s+1}\})$.
Since $u_{j}$ is an endvertex of $Q_{q}$, there exists a subpath $Q$ of $Q_{q}$ such that $u_{j}$ is an endvertex of $Q$ and $|V(Q)|=n$.
Since $|S|\geq 2$, we can take a vertex $v\in S$ as follows:
If $q\neq s+1$, let $v\in S\setminus \{u_{j}\}$; otherwise (i.e., $j=j^{*}$), let $v=u_{j_{1}}$.
Then by the definition of $Q_{p}~(1\leq p\leq s+1)$, $N_{G}(v)\cap V(Q_{q})=\emptyset $.
Since $2n_{0}\geq 2n\geq n+1$, the vertices $x_{i}$ with $2\leq i\leq n+1$ have been defined, and hence this implies that $\{x_{1},v,x_{2},x_{3},\ldots ,x_{n+1}\}\cup V(Q)$ induces a copy of $F^{(1)}_{n,n}$ in $G$, which is a contradiction.
\qed

\begin{lem}
\label{lem2.3}
Let $i$ be an integer with $n_{0}+1\leq i\leq m-n_{0}$, and let $y\in Y_{i}$.
Then the following hold:
\begin{enumerate}
\item[{\upshape(i)}]
If $yu_{i+1}\notin E(G)$, then $N_{G}(y)\cap V(P)=\{u_{i},u_{i+2}\}$.
\item[{\upshape(ii)}]
We have $Y_{m-n_{0}}=\emptyset $.
\item[{\upshape(iii)}]
If $G$ is $F^{(3)}_{n,n}$-free, then $yu_{i+1}\in E(G)$.
\end{enumerate} 
\end{lem}
\proof
\begin{enumerate}
\item[{\upshape(i)}]
Suppose that $yu_{i+1}\notin E(G)$ and $N_{G}(y)\cap V(P)\neq \{u_{i},u_{i+2}\}$.
Let $k=\max\{j:n_{0}+1\leq j\leq m-n_{0},~yu_{j}\in E(G)\}$.
If $k=i$ (i.e., $N_{G}(y)\cap V(P)=\{u_{i}\}$), then
$$
\{u_{i},y,u_{i-1},u_{i-2},\ldots ,u_{i-n},u_{i+1},u_{i+2},\ldots ,u_{i+n}\}
$$
induces a copy of $F^{(1)}_{n,n}$ in $G$, which is a contradiction.
Since $yu_{i+1}\notin E(G)$ and $N_{G}(y)\cap V(P)\neq \{u_{i},u_{i+2}\}$, this forces $k\geq i+3$.
Then
$$
\{u_{i},u_{i+1},u_{i-1},u_{i-2},\ldots ,u_{i-n},y,u_{k},u_{k+1},\ldots ,u_{k+n-2}\}
$$
induces a copy of $F^{(1)}_{n,n}$ in $G$, which is a contradiction.

\item[{\upshape(ii)}]
By (i), if there exists a vertex $y\in Y_{m-n_{0}}$, then it follows that $yu_{m-n_{0}+1}\in E(G)$ or $N_{G}(y)\cap V(P)=\{u_{m-n_{0}},u_{m-n_{0}+2}\}$, and in particular, $N_{G}(y)\cap X_{0}\neq \emptyset $, which contradicts the definition of $Y$.
Thus we have $Y_{m-n_{0}}=\emptyset $.

\item[{\upshape(iii)}]
Suppose that $G$ is $F^{(3)}_{n,n}$-free and $yu_{i+1}\notin E(G)$.
Then it follows from (i) that $N_{G}(y)\cap V(P)=\{u_{i},u_{i+2}\}$, and hence
$$
\{y,u_{i+1},u_{i},u_{i-1},\ldots ,u_{i-n+1},u_{i+2},u_{i+3},\ldots ,u_{i+n+1}\}
$$
induces a copy of $F^{(3)}_{n,n}$ in $G$, which is a contradiction.
\qed
\end{enumerate}

\begin{lem}
\label{lem2.4}
We have $V(G)=V(P)\cup N_{G}(V(P))\cup (\bigcup _{2\leq i\leq 2n_{0}-1}X_{i})$.
\end{lem}
\proof
Suppose that $V(G)\neq V(P)\cup N_{G}(V(P))\cup (\bigcup _{2\leq i\leq 2n_{0}-1}X_{i})$.
Since $G$ is connected, there exists a vertex $z\in V(G)\setminus (V(P)\cup N_{G}(V(P))\cup (\bigcup _{2\leq i\leq 2n_{0}-1}X_{i}))$ adjacent to a vertex $y\in V(P)\cup N_{G}(V(P))\cup (\bigcup _{2\leq i\leq 2n_{0}-1}X_{i})$ in $G$.
By Lemma~\ref{lem2.2} and the definition of $X_{i}$ and $Y$, this implies that $y\in Y$.
Let $i$ be the integer such that $y\in Y_{i}$.
Then by Lemma~\ref{lem2.3}(ii), $n_{0}+1\leq i\leq m-n_{0}-1$.
Let $k=\max\{j:n_{0}+1\leq j\leq m-n_{0},~yu_{j}\in E(G)\}$.
By Lemma~\ref{lem2.3}(i), $k\geq i+1$.
If $k=i+1$, then
$$
\{y,z,u_{i},u_{i-1},\ldots ,u_{i-n+1},u_{i+1},u_{i+2},\ldots ,u_{i+n}\}
$$
induces a copy of $F^{(2)}_{n,n}$ in $G$; if $k\geq i+2$, then
$$
\{y,z,u_{i},u_{i-1},\ldots ,u_{i-n+1},u_{k},u_{k+1},\ldots ,u_{k+n-1}\}
$$
induces a copy of $F^{(1)}_{n,n}$ in $G$.
In either case, we obtain a contradiction.
\qed

Now we recursively define the values $\alpha _{i}~(i\geq 0)$ as follows:
Let $\alpha _{0}=2\lceil \frac{n_{0}}{2} \rceil $, and for $i$ with $i\geq 1$, let $\alpha _{i}=(n-1)R(n,\alpha _{i-1}+1)-1$.

\begin{lem}
\label{lem-alpha}
For an integer $i$ with $i\geq 0$, $\alpha (G[X_{i}])\leq \alpha _{i}$.
\end{lem}
\proof
We proceed by induction on $i$.
If $|V(P)|\leq 2n_{0}$, then $G[X_{0}]$ equals to $P$, and hence $\alpha (G[X_{0}])=\alpha (P)=\lceil \frac{|V(P)|}{2} \rceil \leq n_{0}\leq \alpha _{0}$; if $|V(P)|\geq 2n_{0}+1$, then $G[X_{0}]$ consists of two components each of which is a path of order $n_{0}$, and hence $\alpha (G[X_{0}])=2\lceil \frac{n_{0}}{2} \rceil =\alpha _{0}$.
In either case, we have $\alpha (G[X_{0}])\leq \alpha _{0}$.
Thus we may assume that $i\geq 1$, and suppose that $\alpha (G[X_{i-1}])\leq \alpha _{i-1}$.
Since $X_{i}\subseteq N_{G}(X_{i-1})$, it follows from Lemma~\ref{lem2.1} that $\alpha (G[X_{i}])\leq \alpha (G[N_{G}(X_{i-1})])\leq (n-1)R(n,\alpha _{i-1}+1)-1=\alpha _{i}$, as desired.
\qed

Note that the value $\sum _{1\leq i\leq 2n_{0}-1}\alpha _{i}$ is a constant depending on $n$ only.
Thus, considering Lemmas~\ref{lem2.0}, \ref{lem2.4} and \ref{lem-alpha}, it suffices to show that
\begin{enumerate}
\item[{$\bullet $}]
${\rm pc}(G[V(P)\cup Y])$ is bounded by a constant depending on $n$ only, and
\item[{$\bullet $}]
if $G$ is $\{F^{(3)}_{n,n},F^{(4)}_{n,n}\}$-free, then ${\rm pp}(G[V(P)\cup Y])$ is bounded by a constant depending on $n$ only.
\end{enumerate}
Hence the following lemma completes the proof of Theorem~\ref{thm1}.

\begin{lem}
\label{lem-complete}
\begin{enumerate}
\item[{\upshape(i)}]
We have ${\rm pc}(G[V(P)\cup Y])\leq \max\{3n-6,1\}$.
\item[{\upshape(ii)}]
If $G$ is $\{F^{(3)}_{n,n},F^{(4)}_{n,n}\}$-free, then there exists a Hamiltonian path of $G[V(P)\cup Y]$, i.e., ${\rm pp}(G[V(P)\cup Y])=1$.
\end{enumerate}
\end{lem}
\proof
If $Y=\emptyset $, then $P$ is a Hamiltonian path of $G[V(P)\cup Y]$, and hence ${\rm pc}(G[V(P)\cup Y])={\rm pp}(G[V(P)\cup Y])=1$.
Thus we may assume that $Y\neq \emptyset $.
By Lemma~\ref{lem2.3}(ii), $Y_{m-n_{0}}=\emptyset $.

We first prove (i).
Fix an integer $i$ with $n_{0}+1\leq i\leq m-n_{0}-1$.
Let $Y_{i,1}=\{y\in Y_{i}:yu_{i+1}\in E(G)\}$ and $Y_{i,2}=Y_{i}\setminus Y_{i,1}$.
Then by Lemma~\ref{lem2.3}(i), $N_{G}(y)\cap V(P)=\{u_{i},u_{i+2}\}$ for all $y\in Y_{i,2}$.
Let $j\in \{1,2\}$.
If there exists an independent set $U\subseteq Y_{i,j}$ of $G$ with $|U|=n-1$, then $\{u_{i-1},u_{i}\}\cup U$ induces a copy of $K_{1,n}$ in $G$, which is a contradiction.
Thus $\alpha (G[Y_{i,j}])\leq n-2$.
Since $Y\neq \emptyset $, i.e., $Y_{p,q}\neq \emptyset $ for some $p$ and $q$ with $n_{0}+1\leq p\leq m-n_{0}-1$ and $q\in \{1,2\}$, this implies that $n\geq 3$.
By Lemma~\ref{lem2.0}, there exists a path partition $\PP_{i,j}=\{Q^{(1)}_{i,j},Q^{(2)}_{i,j},\ldots ,Q^{(s_{i,j})}_{i,j}\}$ of $G[Y_{i,j}]$ with $s_{i,j}\leq n-2$, where $\PP_{i,j}=\emptyset $ and $s_{i,j}=0$ if $Y_{i,j}=\emptyset $.
For an integer $t$ with $1\leq t\leq n-2$, if $t\leq s_{i,j}$, let $R^{(t)}_{i,j}$ be the path $u_{i}vQ^{(t)}_{i,j}wu_{i+j}$, where $\{v,w\}$ is the set of endvertices of $Q^{(t)}_{i,j}$; otherwise, let $R^{(t)}_{i,j}$ be the path between $u_{i}$ and $u_{i+j}$ on $P$ (i.e., $R^{(t)}_{i,1}=u_{i}u_{i+1}$ and $R^{(t)}_{i,2}=u_{i}u_{i+1}u_{i+2}$).
We define the value $\xi _{2}$ (resp. $\xi _{3}$) with $\xi _{2}=m-n_{0}$ or $\xi _{2}=m-n_{0}-1$ (resp. $\xi _{3}=m-n_{0}-1$ or $\xi _{3}=m-n_{0}$) according as $m$ is odd or even.
Let
\begin{align*}
R^{(t)}_{1} &= u_{1}u_{2}\cdots u_{n_{0}+1}R^{(t)}_{n_{0}+1,1}u_{n_{0}+2}R^{(t)}_{n_{0}+2,1}u_{n_{0}+3}\cdots u_{m-n_{0}-1}R^{(t)}_{m-n_{0}-1,1}u_{m-n_{0}}u_{m-n_{0}+1}\cdots u_{m},\\
R^{(t)}_{2} &= u_{1}u_{2}\cdots u_{n_{0}+1}R^{(t)}_{n_{0}+1,2}u_{n_{0}+3}R^{(t)}_{n_{0}+3,2}u_{n_{0}+5}\cdots u_{\xi _{2}-2}R^{(t)}_{\xi _{2}-2,2}u_{\xi _{2}}u_{\xi _{2}+1}\cdots u_{m},\mbox{ and}\\
R^{(t)}_{3} &= u_{1}u_{2}\cdots u_{n_{0}+2}R^{(t)}_{n_{0}+2,2}u_{n_{0}+4}R^{(t)}_{n_{0}+4,2}u_{n_{0}+6}\cdots u_{\xi _{3}-2}R^{(t)}_{\xi _{3}-2,2}u_{\xi _{3}}u_{\xi _{3}+1}\cdots u_{m}.
\end{align*}
Then we easily verify that $\{R^{(t)}_{a}:a\in \{1,2,3\},~1\leq t\leq n-2\}$ is a path cover of $G[V(P)\cup Y]$ having cardinality at most $3(n-2)$, which proves (i).

Next we prove (ii).
Suppose that $G$ is $\{F^{(3)}_{n,n},F^{(4)}_{n,n}\}$-free.
We start with the following claim.

\begin{claim}
\label{cl-comp-1}
For an integer $i$ with $n_{0}+1\leq i\leq m-n_{0}-1$, $\{u_{i},u_{i+1}\}\cup Y_{i}$ is a clique of $G$.
\end{claim}
\proof
Suppose that there exist two vertices $y,y'\in \{u_{i},u_{i+1}\}\cup Y_{i}$ with $yy'\notin E(G)$.
By the definition of $Y_{i}$ and Lemma~\ref{lem2.3}(iii), every vertex in $Y_{i}$ is adjacent to both $u_{i}$ and $u_{i+1}$ in $G$.
Thus $y,y'\in Y_{i}$.
Recall that $N_{G}(Y)\cap V(P)\subseteq \{u_{j}:n_{0}+1\leq j\leq m-n_{0}\}$.
Let $k=\max\{j:n_{0}+1\leq j\leq m-n_{0},~N_{G}(u_{j})\cap \{y,y'\}\neq \emptyset \}$.
We may assume that $yu_{k}\in E(G)$.
Note that $k\geq i+1$.
If $k=i+1$, then
$$
\{y,y',u_{i},u_{i-1},\ldots ,u_{i-n+1},u_{i+1},u_{i+2},\ldots ,u_{i+n}\}
$$
induces a copy of $F^{(4)}_{n,n}$ in $G$, which is a contradiction.
Thus $k\geq i+2$.
If $y'u_{k}\in E(G)$, then
$$
\{y,y',u_{i},u_{i-1},\ldots ,u_{i-n+1},u_{k},u_{k+1},\ldots ,u_{k+n-1}\}
$$
induces a copy of $F^{(3)}_{n,n}$ in $G$; if $y'u_{k}\notin E(G)$, then
$$
\{u_{i},y',u_{i-1},u_{i-2},\ldots ,u_{i-n},y,u_{k},u_{k+1},\ldots ,u_{k+n-2}\}
$$
induces a copy of $F^{(1)}_{n,n}$ in $G$.
In either case, we obtain a contradiction.
\qed

For an integer $i$ with $n_{0}+1\leq i\leq m-n_{0}-1$, it follows from Claim~\ref{cl-comp-1} that there exists a Hamiltonian path $R_{i}$ of $G[\{u_{i},u_{i+1}\}\cup Y_{i}]$ with the endvertices $u_{i}$ and $u_{i+1}$.
Then
$$
u_{1}u_{2}\cdots u_{n_{0}+1}R_{n_{0}+1}u_{n_{0}+2}R_{n_{0}+2}u_{n_{0}+3}\cdots u_{m-n_{0}-1}R_{m-n_{0}-1}u_{m-n_{0}}u_{m-n_{0}+1}\cdots u_{m}
$$
is a Hamiltonian path of $G[V(P)\cup (\bigcup _{n_{0}+1\leq i\leq m-n_{0}-1}Y_{i})]~(=G[V(P)\cup Y])$, as desired.
\qed

\subsection{The ``only if'' parts of Theorem~\ref{mainthm}}\label{sec2.2}

Let $s\geq 2$ and $t\geq 3$ be integers, and let $Q_{i}=u^{(1)}_{i}u^{(2)}_{i}\cdots u^{(t)}_{i}~(1\leq i\leq s)$ be $s$ pairwise vertex-disjoint paths.
We define four graphs.
\begin{enumerate}[{$\bullet $}]
\item
Let $H^{(1)}_{s,t}$ be the graph obtained from the union of the paths $Q_{1},\ldots ,Q_{s}$ by adding $2(s-1)$ vertices $v_{i},w_{i}~(1\leq i\leq s-1)$ and $3(s-1)$ edges $v_{i}w_{i},v_{i}u^{(t)}_{i},v_{i}u^{(1)}_{i+1}~(1\leq i\leq s-1)$.
\item
Let $H^{(2)}_{s,t}$ be the graph obtained from $H^{(1)}_{s,t}$ by adding $s-1$ edges $u^{(t)}_{i}u^{(1)}_{i+1}~(1\leq i\leq s-1)$.
\item
Let $H^{(3)}_{s,t}$ be the graph obtained from the union of the paths $Q_{1},\ldots ,Q_{s}$ by adding $2(s-1)$ vertices $v_{i},w_{i}~(1\leq i\leq s-1)$ and $4(s-1)$ edges $v_{i}u^{(t)}_{i},v_{i}u^{(1)}_{i+1},w_{i}u^{(t)}_{i},w_{i}u^{(1)}_{i+1}~(1\leq i\leq s-1)$.
\item
Let $H^{(4)}_{s,t}$ be the graph obtained from $H^{(3)}_{s,t}$ by adding $s-1$ edges $u^{(t)}_{i}u^{(1)}_{i+1}~(1\leq i\leq s-1)$.
\end{enumerate}

\begin{lem}
\label{lem-H12}
We have ${\rm pc}(H^{(1)}_{s,t})={\rm pc}(H^{(2)}_{s,t})=\lceil \frac{s+1}{2} \rceil$.
\end{lem}
\proof
Note that $u^{(1)}_{1},u^{(t)}_{s},w_{i}~(1\leq i\leq s-1)$ have degree one in $H^{(2)}_{s,t}$.
Since a path contains at most two vertices of degree at most one, ${\rm pc}(G)\geq \lceil \frac{l}{2} \rceil $ for every graph $G$ where $l$ is the number of the vertices of $G$ having degree one.
In particular, we have
\begin{align}
{\rm pc}(H^{(2)}_{s,t})\geq \left\lceil \frac{s+1}{2} \right\rceil .\label{lem-H12-1}
\end{align}
If $s$ is odd, let
$$
\PP=\left\{H^{(1)}_{s,t}-\{w_{j}:1\leq j\leq s-1\},~w_{2i-1}v_{2i-1}u^{(1)}_{2i}Q_{2i}u^{(t)}_{2i}v_{2i}w_{2i}:1\leq i\leq \frac{s-1}{2}\right\};
$$
if $s$ is even, let
$$
\PP=\left\{H^{(1)}_{s,t}-\{w_{j}:1\leq j\leq s-1\},~H^{(1)}_{s,t}[\{w_{s-1}\}],~w_{2i-1}v_{2i-1}u^{(1)}_{2i}Q_{2i}u^{(t)}_{2i}v_{2i}w_{2i}:1\leq i\leq \frac{s-2}{2}\right\}.
$$
Then we verify that $\PP$ is a path cover of $H^{(1)}_{s,t}$ with $|\PP|=\lceil \frac{s+1}{2} \rceil $.
Furthermore, since $H^{(1)}_{s,t}$ is a spanning subgraph of $H^{(2)}_{s,t}$, a path cover of $H^{(1)}_{s,t}$ is also a path cover of $H^{(2)}_{s,t}$, and hence ${\rm pc}(H^{(2)}_{s,t})\leq {\rm pc}(H^{(1)}_{s,t})\leq \lceil \frac{s+1}{2} \rceil $.
This together with (\ref{lem-H12-1}) leads to the desired conclusion.
\qed

\begin{lem}
\label{lem-H34}
We have ${\rm pp}(H^{(3)}_{s,t})={\rm pp}(H^{(4)}_{s,t})=s$.
\end{lem}
\proof
We first prove that
\begin{align}
{\rm pp}(H^{(4)}_{s,t})\geq s.\label{lem-H34-1}
\end{align}
Let $\PP$ be a path partition of $H^{(4)}_{s,t}$.
It suffices to show that $|\PP|\geq s$.
For each $i$ with $1\leq i\leq s$, let $R_{i}$ be the unique element of $\PP$ containing $u^{(2)}_{i}$.
We remark that $R_{i}$ might equal to $R_{j}$ for some $1\leq i<j\leq s$.
Let $I=\{i:1\leq i\leq s-1,~R_{i}=R_{i+1}\}$, and write $I=\{i_{1},i_{2},\ldots ,i_{h}\}$ with $i_{1}<i_{2}<\ldots <i_{h}$ where $h=0$ if $I=\emptyset $.
For integers $i$ and $i'$ with $1\leq i<i'\leq s$, any paths of $H^{(4)}_{s,t}$ joining $u^{(2)}_{i}$ and $u^{(2)}_{i'}$ contain every vertex in $\{u^{(2)}_{j}:i<j<i'\}$.
This implies that if $R_{i}=R_{i'}$ with $1\leq i<i'\leq s$, then $i'-i+1$ paths $R_{j}~(i\leq j\leq i')$ are equal.
In particular, we have $|\{R_{i}:1\leq i\leq s\}|=s-h$.

Fix an integer $l$ with $1\leq l\leq h$.
Then for every path $R$ of $H^{(4)}_{s,t}$ joining $u^{(2)}_{i_{l}}$ and $u^{(2)}_{i_{l}+1}$, we easily verify that
\begin{enumerate}
\item[$\bullet $]
$\{u^{(t)}_{i_{l}},u^{(1)}_{i_{l}+1}\}\subseteq V(R)$, and
\item[$\bullet $]
$v_{i_{l}}\notin V(R)$ or $w_{i_{l}}\notin V(R)$.
\end{enumerate}
Since $v_{i_{l}}w_{i_{l}}\notin E(H^{(4)}_{s,t})$, this implies that there exists an element $R'_{i_{l}}$ of $\PP$ such that either $V(R'_{i_{l}})=\{v_{i_{l}}\}$ or $V(R'_{i_{l}})=\{w_{i_{l}}\}$.
Therefore
\begin{align*}
|\PP| &\geq |\{R_{i}:1\leq i\leq s\}\cup \{R'_{i_{j}}:1\leq j\leq h\}|\\
&= |\{R_{i}:1\leq i\leq s\}|+|\{R'_{i_{j}}:1\leq j\leq h\}|\\
&= (s-h)+h\\
&=s,
\end{align*}
which proves (\ref{lem-H34-1}).

Since
$$
\PP'=\{H^{(3)}_{s,t}-\{w_{j}:1\leq j\leq s-1\},~H^{(3)}_{s,t}[\{w_{i}\}]:1\leq i\leq s-1\}
$$
is a path partition of $H^{(3)}_{s,t}$ with $|\PP'|=s$.
Furthermore, since $H^{(3)}_{s,t}$ is a spanning subgraph of $H^{(4)}_{s,t}$, a path partition of $H^{(3)}_{s,t}$ is also a path partition of $H^{(4)}_{s,t}$, and hence ${\rm pp}(H^{(4)}_{s,t})\leq {\rm pp}(H^{(3)}_{s,t})\leq s$.
This together with (\ref{lem-H34-1}) leads to the desired conclusion.
\qed

Now we prove the following proposition, which gives the ``only if'' parts of Theorem~\ref{mainthm}.

\begin{prop}
\label{prop-onlyif}
Let $\HH$ be a finite family of connected graphs.
\begin{enumerate}
\item[{\upshape(i)}]
If $\HH$ satisfies (A1), then $\HH\leq \{K_{1,n},K^{*}_{n},F^{(1)}_{n,n},F^{(2)}_{n,n}\}$ for an integer $n\geq 2$.
\item[{\upshape(ii)}]
If $\HH$ satisfies (A2), then $\HH\leq \{K_{1,n},K^{*}_{n},F^{(1)}_{n,n},F^{(2)}_{n,n},F^{(3)}_{n,n},F^{(4)}_{n,n}\}$ for an integer $n\geq 2$.
\end{enumerate}
\end{prop}
\proof
Since $\HH$ is a finite family, the value $p=\max\{|V(H)|:H\in \HH\}$ is well-defined.
If $p\leq 2$, then the desired conclusions trivially hold.
Thus we may assume that $p\geq 3$.

We first suppose that $\HH$ satisfies (A1), and show that (i) holds.
There exists a constant $c_{1}=c_{1}(\HH)$ such that ${\rm pc}(G)\leq c_{1}$ for every connected $\HH$-free graph $G$.
Since ${\rm pc}(K_{1,2c_{1}+1})=c_{1}+1$ and ${\rm pc}(K^{*}_{2c_{1}+1})=c_{1}+1$, neither $K_{1,2c_{1}+1}$ nor $K^{*}_{2c_{1}+1}$ is $\HH$-free.
This implies that
\begin{align}
\HH\leq \{K_{1,2c_{1}+1},K^{*}_{2c_{1}+1}\}.\label{prop-onlyif-(i)1}
\end{align}
For each $i\in \{1,2\}$, it follows from Lemma~\ref{lem-H12} that ${\rm pc}(H^{(i)}_{2c_{1},p})=\lceil \frac{2c_{1}+1}{2} \rceil =c_{1}+1$, and hence $H^{(i)}_{2c_{1},p}$ is not $\HH$-free, i.e., $H^{(i)}_{2c_{1},p}$ contains an induced subgraph $A_{i}$ isomorphic to an element of $\HH$.
Since $|V(A_{i})|\leq p$, we have
\begin{enumerate}
\item[$\bullet $]
$|\{j:1\leq j\leq 2c_{1},~V(A_{i})\cap V(Q_{j})\neq \emptyset \}|\leq 2$, and
\item[$\bullet $]
$|\{j:1\leq j\leq 2c_{1}-1,~V(A_{i})\cap \{v_{j},w_{j}\}\neq \emptyset \}|\leq 1$.
\end{enumerate}
This implies that $A_{i}$ is an induced copy of $F^{(i)}_{p,p}$, and hence
\begin{align}
\HH\leq \{F^{(1)}_{p,p},F^{(2)}_{p,p}\}.\label{prop-onlyif-(i)2}
\end{align}
Let $n=\max\{2c_{1}+1,p\}$.
Then by (\ref{prop-onlyif-(i)1}) and (\ref{prop-onlyif-(i)2}), $\HH\leq \{K_{1,n},K^{*}_{n},F^{(1)}_{n,n},F^{(2)}_{n,n}\}$, which proves (i).

Next we suppose that $\HH$ satisfies (A2), and show that (ii) holds.
There exists a constant $c_{2}=c_{2}(\HH)$ such that ${\rm pp}(G)\leq c_{2}$ for every connected $\HH$-free graph $G$.
Since ${\rm pp}(G)\geq {\rm pc}(G)$ for all graphs $G$, $\HH$ also satisfies (A1).
Hence by (i), there exists an integer $m\geq 2$ such that
\begin{align}
\HH\leq \{K_{1,m},K^{*}_{m},F^{(1)}_{m,m},F^{(2)}_{m,m}\}.\label{prop-onlyif-(ii)1}
\end{align}
For each $i\in \{3,4\}$, it follows from Lemma~\ref{lem-H34} that ${\rm pp}(H^{(i)}_{c_{2}+1,p})=c_{2}+1$, and hence $H^{(i)}_{c_{2}+1,p}$ is not $\HH$-free, i.e., $H^{(i)}_{c_{2}+1,p}$ contains an induced subgraph $B_{i}$ isomorphic to an element of $\HH$.
Since $|V(B_{i})|\leq p$, we have
\begin{enumerate}
\item[$\bullet $]
$|\{j:1\leq j\leq c_{2}+1,~V(B_{i})\cap V(Q_{j})\neq \emptyset \}|\leq 2$, and
\item[$\bullet $]
$|\{j:1\leq j\leq c_{2},~V(B_{i})\cap \{v_{j},w_{j}\}\neq \emptyset \}|\leq 1$.
\end{enumerate}
This implies that $B_{i}$ is an induced copy of $F^{(i)}_{p,p}$, and hence
\begin{align}
\HH\leq \{F^{(3)}_{p,p},F^{(4)}_{p,p}\}.\label{prop-onlyif-(ii)2}
\end{align}
Let $n'=\max\{m,p\}$.
Then by (\ref{prop-onlyif-(ii)1}) and (\ref{prop-onlyif-(ii)2}), $\HH\leq \{K_{1,n'},K^{*}_{n'},F^{(1)}_{n',n'},F^{(2)}_{n',n'},F^{(3)}_{n',n'},F^{(4)}_{n',n'}\}$, which proves (ii).
\qed

\section{Proof of Theorem~\ref{mainthm2}}\label{sec3}

In this section, we prove Theorem~\ref{mainthm2}.
We start with the following lemma, which is an analogy of Lemma~\ref{lem2.0}.

\begin{lem}
\label{lem3.0}
For a graph $G$, ${\rm cp}(G)\leq R(\alpha (G)+1,\alpha (G)+1)-1$.
\end{lem}
\proof
Let $\PP$ be a cycle partition of $G$ with $|\PP|={\rm cp}(G)$, and write $\PP=\{Q_{i}:1\leq i\leq {\rm cp}(G)\}$.
By way of contradiction, suppose that ${\rm cp}(G)\geq R(\alpha (G)+1,\alpha (G)+1)$.
For each integer $i$ with $1\leq i\leq {\rm cp}(G)$, we define vertices $x_{i}$ and $y_{i}$ of $Q_{i}$ as follows:
If either $Q_{i}\simeq K_{2}$ or $Q_{i}$ is a cycle, let $x_{i}$ and $y_{i}$ be vertices of $Q_{i}$ with $x_{i}y_{i}\in E(Q_{i})$; if $Q_{i}\simeq K_{1}$, let $x_{i}=y_{i}=u$ where $u$ is the unique vertex of $Q_{i}$.
For integers $i$ and $j$ with $1\leq i<j\leq {\rm cp}(G)$, if $\{x_{i}x_{j},y_{i}y_{j}\}\subseteq E(G)$, then we easily verify that there exists a spanning subgraph $Q$ of $G[V(Q_{i})\cup V(Q_{j})]$ such that either $Q\simeq K_{2}$ or $Q$ is a cycle, and hence $\PP'=(\PP\setminus \{Q_{i},Q_{j}\})\cup \{Q\}$ is a cycle partition of $G$ with $|\PP'|={\rm cp}(G)-1$, which contradicts the definition of the cycle partition number.
Thus if $x_{i}x_{j}\in E(G)$, then $y_{i}y_{j}\notin E(G)$.

Let $K$ be the complete graph on $\{1,2,\ldots ,{\rm cp}(G)\}$, and color all edges of $K$ by red or blue as follows:
For integers $i$ and $j$ with $1\leq i<j\leq {\rm cp}(G)$, if $x_{i}x_{j}\notin E(G)$, we color the edge $ij$ of $K$ by red; if $x_{i}x_{j}\in E(G)$ and $y_{i}y_{j}\notin E(G)$, we color the edge $ij$ of $K$ by blue.
Since $|V(K)|={\rm cp}(G)\geq R(\alpha (G)+1,\alpha (G)+1)$, there exists a monochromatic clique $I$ of $K$ with $|I|=\alpha (G)+1$.
If $I$ is a red clique of $K$, then $\{x_{i}:i\in I\}$ is an independent set of $G$; if $I$ is a blue clique of $K$, then $\{y_{i}:i\in I\}$ is an independent set of $G$.
In either case, we obtain a contradiction.
\qed

The following lemma was implicitly proved in \cite{CFKP}.
(To keep the paper self-contained, we give its proof.)

\begin{lem}[Choi et al.~\cite{CFKP}]
\label{lem3.1}
Let $n\geq 2$ be an integer.
There exists a constant $c=c(n)$ depending on $n$ only such that $\alpha (G)\leq c$ for every connected $\{K_{1,n},K^{*}_{n},P_{n}\}$-free graph $G$.
\end{lem}
\proof
Let $x$ be a vertex of $G$, and for an integer $i$ with $i\geq 0$, let $X_{i}$ be the set of vertices $y$ of $G$ such that the distance between $x$ and $y$ in $G$ is exactly $i$.
Note that $X_{0}=\{x\}$ and $X_{1}=N_{G}(x)$.
Since $G$ is $P_{n}$-free, $X_{i}=\emptyset $ for all $i\geq n-1$.
Since $G$ is connected, this implies that
\begin{align}
V(G)=\bigcup _{0\leq i\leq n-2}X_{i}.\label{cond-lem3.1-0}
\end{align}
We recursively define the values $\alpha _{i}~(i\geq 0)$ as follows:
Let $\alpha _{0}=1$, and for $i$ with $i\geq 1$, let $\alpha _{i}=(n-1)R(n,\alpha _{i-1}+1)-1$.

We prove that
\begin{align}
\mbox{$\alpha (G[X_{i}])\leq \alpha _{i}$ for an integer $i$ with $0\leq i\leq n-2$.}\label{cond-lem3.1-1}
\end{align}
We proceed by induction on $i$.
Since $\alpha (G[X_{0}])=1=\alpha _{0}$, we may assume that $i\geq 1$ and $\alpha (G[X_{i-1}])\leq \alpha _{i-1}$.
Since $X_{i}\subseteq N_{G}(X_{i-1})$, it follows from Lemma~\ref{lem2.1} that $\alpha (G[X_{i}])\leq \alpha (G[N_{G}(X_{i-1})])\leq (n-1)R(n,\alpha _{i-1}+1)-1=\alpha _{i}$, as desired.

By (\ref{cond-lem3.1-0}) and (\ref{cond-lem3.1-1}), we have $\alpha (G)\leq \sum _{0\leq i\leq n-2}\alpha (G[X_{i}])\leq \sum _{0\leq i\leq n-2}\alpha _{i}$.
Since the value $\sum _{0\leq i\leq n-2}\alpha _{i}$ is a constant depending on $n$ only, we obtain the desired conclusion.
\qed

\medbreak\noindent\textit{Proof of Theorem~\ref{mainthm2}.}\quad
By the definition of cycle cover and cycle partition, ``(ii) $\Longrightarrow $ (i)'' clearly holds.

We show that ``(iii) $\Longrightarrow $ (ii)'' holds.
Let $n\geq 2$ be an integer, and let $c=c(n)$ be the constant as in Lemma~\ref{lem3.1}.
It suffices to show that there exists a constant $c_{1}=c_{1}(n)$ depending on $n$ only such that ${\rm cp}(G)\leq c_{1}$ for every connected $\{K_{1,n},K^{*}_{n},P_{n}\}$-free graph $G$.
By the definition of $c(n)$, we have $\alpha (G)\leq c$.
This together with Lemma~\ref{lem3.0} leads to ${\rm cp}(G)\leq R(\alpha (G)+1,\alpha (G)+1)-1\leq R(c+1,c+1)-1$.
Since $R(c+1,c+1)-1$ is a constant depending on $n$ only, we obtain the desired conclusion.

Finally, we show that ``(i) $\Longrightarrow $ (iii)'' holds, which completes the proof of Theorem~\ref{mainthm2}.
Suppose that a family $\HH$ of connected graphs satisfies (A'1).
Then there exists a constant $c_{1}=c_{1}(\HH)$ such that ${\rm cc}(G)\leq c_{1}$ for every connected $\HH$-free graph $G$.
Since ${\rm cc}(K_{1,c_{1}+1})=c_{1}+1$, ${\rm cc}(K^{*}_{c_{1}+1})=c_{1}+1$ and ${\rm cc}(P_{2c_{1}+1})=\lceil \frac{2c_{1}+1}{2} \rceil =c_{1}+1$, none of $K_{1,c_{1}+1}$, $K^{*}_{c_{1}+1}$ and $P_{2c_{1}+1}$ is $\HH$-free.
This implies that $\HH\leq \{K_{1,c_{1}+1},K^{*}_{c_{1}+1},P_{2c_{1}+1}\}$, and hence $\HH\leq \{K_{1,2c_{1}+1},K^{*}_{2c_{1}+1},P_{2c_{1}+1}\}$, which leads (iii).
\qed

\section*{Acknowledgment}

This work was partially supported by JSPS KAKENHI Grant number JP20K03720 (to S.C) and JSPS KAKENHI Grant number JP18K13449 (to M.F).


\begin{thebibliography}{99}

\bibitem{CFKP}
I.~Choi, M.~Furuya, R.~Kim and B.~Park,
A Ramsey-type theorem for the matching number regarding connected graphs,
Discrete Math. {\bf 343} (2020), 111648.

\bibitem{D}
R.~Diestel,
``Graph Theory'' (5th edition), Graduate Texts in Mathematics \textbf{173},
Springer, Berlin (2017).

\bibitem{DGJ}
D.~Duffus, R.J.~Gould and M.S.~Jacobson,
Forbidden subgraphs and the hamiltonian theme,
The Theory and Applications of Graphs, pp. 297--316, Wiley, New York (1981).

\bibitem{EL}
H.~Enomoto and H.~Li,
Partition of a graph into cycles and degenerated cycles,
Discrete Math. {\bf 276} (2004), 177--181.

\bibitem{FG}
R.J.~Faudree and R.J.~Gould,
Characterizing forbidden pairs for Hamiltonian properties,
Discrete Math. {\bf 173} (1997), 45--60.

\bibitem{Fj}
S.~Fujita,
Degree conditions for the partition of a graph into cycles, edges and isolated vertices,
Discrete Math. {\bf 309} (2009), 3534--3540.

\bibitem{FKLOPS}
S.~Fujita, K.~Kawarabayashi, C.~L.~Lucchesi, K.~Ota, M.~Plummer and A.~Saito,
A pair of forbidden subgraphs and perfect matchings,
J. Combin. Theory Ser. B {\bf 96} (2006), 315--324.

\bibitem{F}
M.~Furuya,
Forbidden subgraphs for constant domination number,
Discrete Math. Theor. Comput. Sci. {\bf 20(1)} (2018), Paper No. 19.

\bibitem{GM}
T.~Gallai and A.N.~Milgram,
Verallgemeinerung eines graphentheoretischen Satzes von Redei,
Acta Sci. Math. {\bf 21} (1960), 181--186.

\bibitem{GH1}
R.J.~Gould and J.M.~Harris,
Forbidden triples of subgraphs and traceability,
Congr. Numer. {\bf 108} (1995), 183--192.

\bibitem{GH2}
R.J.~Gould and J.M.~Harris,
Traceability in graphs with forbidden triples of subgraphs,
Discrete Math. {\bf 189} (1998), 123--132.

\bibitem{GH3}
R.J.~Gould and J.M.~Harris,
Forbidden triples and traceability: a characterization,
Discrete Math. {\bf 203} (1999), 101--120.

\bibitem{H}
J.~Han,
On vertex-disjoint paths in regular graphs,
Electron. J. Combin. {\bf 25} (2018), Paper No. 2.12.

\bibitem{I}
S.~Ishizuka,
Closure, path-factors and path coverings in claw-free graphs,
Ars Combin. {\bf 50} (1998), 115--128.

\bibitem{MM}
C.~Magnant and D.M.~Martin,
A note on the path cover number of regular graphs,
Australas. J. Combin. {\bf 43} (2009), 211--217.

\bibitem{R}
B.~Reed,
Paths, stars and the number three,
Combin. Probab. Comput. {\bf 5} (1996), 277--295.

\bibitem{Y}
G.~Yu,
Covering $2$-connected $3$-regular graphs with disjoint paths,
J. Graph Theory {\bf 88} (2018), 385--401.




\end{thebibliography}
\end{document}